\documentclass[10pt,a4paper]{amsart}
\usepackage{amssymb,amsmath,amsthm,amsfonts,amstext}
\usepackage{array}
\usepackage{longtable}
\usepackage{textcomp}
\usepackage{multirow}
\usepackage{graphicx}
\usepackage{esint}
\usepackage{color}

\theoremstyle{plain}
\newtheorem{thm}{Theorem}

\newtheorem{cor}[thm]{Corollary}
\theoremstyle{definition}

\begin{document}

\title{modeling extreme values by the residual coefficient of variation}

\author{Joan del Castillo \and Maria Padilla}

 \address{Joan del Castillo. Department of Mathematics, Universitat Aut\`onoma
de Barcelona, 08193, Cerdanyola del Valles (Barcelona), Spain. Corresponding author:} 
\email{castillo@mat.uab.cat}
\address{Maria Padilla. Department of Mathematics. Universitat Aut\`onoma de
Barcelona. 08193, Cerdanyola del Valles (Barcelona). Spain.}
\date{}


\maketitle

\begin{abstract}
The possibilities of the use of the coefficient of variation over a high threshold in tail modelling are discussed. The paper also considers multiple threshold tests for a generalized Pareto distribution, together with a threshold selection algorithm. One of the main contributions is to extend the methodology based on moments to all distributions, even without finite moments. These techniques are applied to Danish fire insurance losses.
\end{abstract}

\section{Introduction}

Fisher \& Tippett~\cite{Fisher28} 
and Gnedenko~\cite{Gnedenko43}  
show that, under regularity
conditions, the limit distribution for the normalized maximum of a
sequence of independent and identically
distributed (iid) random variable (r.v.) is a member of the generalized extreme value~(GEV) distribution with a cumulative
distribution function
\[
H_{\xi}(x)=\exp\{-(1+\xi x)^{-1/\xi}\},\quad (1+\xi x)>0,
\]
where $\xi$ is called  \emph{extreme value index}. This family of continuous distributions contains the Fr\'echet distribution
($\xi>0$), the Weibull distribution ($\xi<0$), and the Gumbell distribution
($\xi=0,$ as a limit case), see~\cite{McNeil05}.  

The Pickands--Balkema--DeHaan Theorem, see~\cite{Embrechts97} 
and \cite{McNeil05},  
initiated a new way of studying
extreme value theory via distributions above a threshold, which use
more information than the maximum data grouped into blocks. This Theorem is a very widely applicable result
that essentially says that the generalized
Pareto distribution~($\mathrm{GPD}$) is the canonical distribution for
modelling excess losses over high thresholds. The cumulative distribution
function of $\mathrm{GPD}(\xi,\psi)$ is
\begin{equation}
F(x)=1-(1+\xi x/\psi)^{-1/\xi},\label{eq:gpd}
\end{equation}
where $\psi>0$ and $\xi$ are scale and shape parameters. For $\xi>0$
the range of $x$ is~$x>0$, in this case the $\mathrm{GPD}$ is simply the usual
Pareto distribution. The limit case $\xi=0$ corresponds to the
exponential distribution. For $\xi<0$ the range of $x$ is $0<x<\psi/\vert \xi\vert $
and $\mathrm{GPD}$ has bounded support. The shape parameter $\xi$ in GPD corresponds to the extreme value index in GEV. The $\mathrm{GPD}$ has mean $\psi/(1-\xi)$ and variance $\psi^{2}/[(1-\xi)^{2}(1-2\xi)]$
provided $\xi<1/2$.

Let $X$ be a continuous non-negative r.v.\ with
distribution function $F(x)$. For any threshold, $t>0$,
the r.v.\ of the conditional distribution of threshold excesses $X-t$
given $X>t$, denoted $X_{t}=(X-t\mid X>t)$, is called
the \emph{residual distribution} of $X$ over $t$. The cumulative
distribution function of~$X_{t}$, $F_{t}(x),$ is given
by
\begin{equation}
1-F_{t}(x)=(1-F(x+t))/(1-F(t)).\label{eq:thresh}
\end{equation}

The quantity $M(t)=E(X_{t})$ is called the \emph{residual mean} and
$V(t)=\operatorname{var}(X_{t})$ the \emph{residual variance}. The \emph{residual
coefficient of variation}  ($\mathrm{CV}$) is given by 
\begin{equation}
\mathrm{CV}(t)\equiv \mathrm{CV}(X_{t})=\sqrt{V(t)}/M(t),\label{eq:CV1}
\end{equation}
like the usual $\mathrm{CV}$, the function $\mathrm{CV}(t)$ is independent of scale,
that is, if $X$ is multiplied by a positive constant, $\mathrm{CV}(t)$ is
invariant.

The residual distribution of a $\mathrm{GPD}$ is again $\mathrm{GPD}$ and for any threshold
$t>0$, the shape parameter~$\xi$ is invariant, in fact 
\begin{equation}
\mathrm{GPD}_{t}(\xi,\psi)=\mathrm{GPD}(\xi,\psi+\xi t).\label{eq:invshape}
\end{equation}

Note that the residual $\mathrm{CV}$ is independent of the
threshold and the scale parameter, since it is given by
\begin{equation}
\mathrm{CV}(t)=c_{\xi}=\sqrt{1/(1-2\xi)}.\label{eq:cvgpd}
\end{equation}

Gupta and Kirmani~\cite{Gupta00}   
show that the residual $\mathrm{CV}$ characterizes the distribution in the univariate as well as the
bivariate case, provided there is a finite second moment. In the case
of $\mathrm{GPD}$, the residual~$\mathrm{CV}$ is constant and is a one to one transformation
of the extreme value index suggesting its use to estimate this index. 

Castillo \emph{et al.}~\cite{Castillo14} suggest a new tool to identify the tail of a distribution
based on the residual $\mathrm{CV}$, henceforth called
CV-plot, as an alternative to the  \emph{mean excess plot} (ME-plot) that is a commonly
used diagnostic tool in risk analysis to justify fitting a $\mathrm{GPD}$, see \cite{Ghosh10},  \cite{Embrechts97} 
and  \cite{Davison90}. Given a sample $\{ x_{k}\} $ of size~$n$ of positive
numbers, we denote the ordered sample~$\{ x_{(k)}\} $,
so that $x_{(1)}\leq x_{(2)}\leq\dotsb\leq x_{(n)}$.
The CV-plot is the $\mathrm{CV}$ of the residual samples, that is, the
function, $cv(t)$ of the $\mathrm{CV}$ of the \emph{threshold
excesses} $(x_{j}-t)$ for the \emph{exceedances} $\{ x_{j}:x_{j}>t\} $,
over the order statistics, $t=x_{(k)}$, given by
\begin{equation}
k\rightarrow cv(x_{(k)})=\frac{sd\{ x_{j}-x_{(k)}\mid x_{j}>x_{(k)}\} }{\text{mean}\{ x_{j}-x_{(k)}\mid x_{j}>x_{(k)}\} },\label{eq:sample cv}
\end{equation}
where, $k$ $(1\leq k\leq n)$ is the size of the sub-sample removed. This tool has been applied to financial and environmental datasets, see  \cite{CastilloSerra14}.

The $\mathrm{CV}$-plot has some advantages over ME-plot: first,
it does not depend on the scale parameter; second, detecting constant
functions is easier than linear functions, since linear functions
are defined by two parameters and the constants by only one. The uncertainty
is essentially reduced from three to one single parameter.

A unconscientious use of some measures of variation can lead to wrong conclusion, see \cite{Albrecher10}.
A serious problem with the residual coefficient of variation  is the fact that the proposed method
only works when the extreme value index is smaller than 0.25. To fix this, some transformations that relate
light-heavy tails are introduced in Section~\ref{section2}.

Section~\ref{section3} extends  some results of Castillo \emph{et al.}~\cite{Castillo14}  from the exponential distribution to all GPD when the extreme value index is below  0.25. Moreover, multiple threshold tests together with a threshold selection algorithm, designed in a way that avoids subjectivity, are also achieved.
In Section~\ref{exa:Danish}, the approach developed in the previous sections is illustrated using the Danish fire insurance dataset, a highly heavy-tailed, infinite-variance model.

\section{Transformations of heavy-light tails}\label{sec:Transformations}
\label{section2}

The transformations introduced to this section make it possible to estimate the extreme value index using methods based on moments in situations where moments are not finite.

A distribution function $F$ is said to be in the maximum domain of
attraction of $H_{\xi}$, written $F\in D(H_{\xi})$, if under appropriate
normalization the block maxima of a iid sequence of r.v.\ with distribution~$F$ converge to $H_{\xi}$. For a r.v.~$X$ with distribution function~$F$ is also written $X\in D(H_{\xi})$. A positive function~$L$
on $(0,\infty)$ slowly varies at $\infty$ if 
\[
\lim_{x\to \infty}\frac{L(tx)}{L(x)}=1,\quad t>0.
\]

Regularly varying functions can be represented by power functions
multiplied by slowly varying functions, i.e.\ $h(x)\in \mathrm{RV}_{\rho}$
if and only if $h(x)=x^{\rho}L(x)$. 

Gnedenko proved, see \cite[Theorems~7.8 and 7.10]{McNeil05}, 
that the maximum domain of attraction of a Fr\'echet distribution,
with shape parameter~$\xi>0$, is characterized in terms of the tail
function, $\overline{F}(x)=1-F(x)$, by 
\[
F\in D(H_{\xi})\Leftrightarrow\overline{F}\in \mathrm{RV}_{-1/\xi}\quad (\xi>0).
\]

Similarly the maximum domain of attraction of a Weibull distribution,
with shape parameter $\xi<0$, is characterized by
\[
F\in D(H_{\xi})\Leftrightarrow\overline{F}(x_{+}-1/x)\in \mathrm{RV}_{1/\xi}\quad (\xi<0),
\]
 where $x_{+}=\sup\{ x:F(x)<1\} $.

The following result of practical importance is embedded in the previous
characterizations, and which to our knowledge has not been used.

\begin{cor}\label{cor:7}
Let $X$ be a r.v.\ with cumulative distribution function~$F$.
\end{cor}

\begin{enumerate}
\item[(1)] If $X\in D(H_{\xi})$ with $\xi>0$, then $X^{*}=-1/X\in D(H_{-\xi})$. 
\item[(2)] If $X\in D(H_{\xi})$ with $\xi<0$, then $X^{*}=x_{+}-1/X\in D(H_{-\xi})$,
where $x_{+}=\sup\{ x:F(x)<1\} $.
\end{enumerate}

\begin{proof}
(1) The cumulative distribution function of $X^{*}$ is $F^{*}(x)=F(-1/x)$
and $x_{+}=\sup\{ x:F^{*}(x)<1\} =0$. By assumption $\overline{F}(x)=x^{-1/\xi}L(x)$
with $L$ slowly varying at $\infty$, hence $\overline{F^{*}}(x_{+}-1/x)=\overline{F}(x)=x^{1/(-\xi)}L(x)$
and $X^{*}\in D(H_{-\xi})$. 

(2) The tail function of $X^{*}$ is now $\overline{F^{*}}(x)=\overline{F}(x_{+}-1/x)=x^{1/\xi}L(x)$.
Hence, $\overline{F^{*}}(x)\in \mathrm{RV}_{1/\xi}$ and $X^{*}\in D(H_{-\xi})$.
\end{proof}

Corollary \ref{cor:7} provides an asymptotic method and is related
to an exact result in the GEV model: $X$ has Fr\'echet distribution
if and only if $-1/X$ has Weibull distribution with the same extreme value index, but with the sign changed. However, the corresponding result
is not true in $\mathrm{GPD}$, as we discuss below.

For a r.v.\ $X$, the Pickands--Balkema--DeHaan Theorem shows that $X\in D(H_{\xi})$
if and only if the limiting behavior of the residual distribution
of $X$ over $t$, $X_{t}$, is like a $\mathrm{GPD}$ with the same parameter~$\xi$, see 
\cite[Theorem~7.20]{McNeil05}. 
Hence, Corollary~\ref{cor:7} can be used in applied methods of threshold exceedances.

\begin{cor}\label{cor:8}
Let $X$ be a r.v.\ such that the limiting behavior of
the residual distribution of $X$ over a threshold is $\mathrm{GPD}$ with parameter~$\xi$, then the limiting behaviour of the residual distribution of $X^{*}=-1/X$ over a threshold is $\mathrm{GPD}$ with parameter~$-\xi$.
\end{cor}

Corollary \ref{cor:8} enables use of methods to determine the extreme value index
for light tails in heavy tailed distributions and vice versa.
For instance ME-plot and $\mathrm{CV}$-plot can be used to determine the extreme value index in really heavy tailed distributions, see Example~\ref{exa:Danish} below. These asymptotic results can be improved on GPD for practical aplications.

The $\mathrm{GPD}(\xi,\psi)$ distributions are standardized so
that all their observations take positive values. The supports of
the distributions are $(0,\sigma)$, where $\sigma=\infty$
for $\xi\geq0$ and $\sigma=\psi/|\xi|$ for $\xi<0$.
The $\mathrm{GPD}$ distributions can be expanded to include a location parameter
by $Y=X+\mu$. The behavior of $X$ near~$\sigma$ is the same as
that of $Y$ near $\sigma+\mu$. The transformation $X^{*}=-1/X$
is also associated with the origin at zero, but can be generalized
to $Y=-1/(X+c)$, provided $c\geq0$, or $c\leq-\sigma$ , in order
for the transformations to remain monotonous increasing on $(0,\sigma)$.
The following result examines these transformations on~$\mathrm{GPD}$.

\begin{thm}\label{thm:9}
Let $X$ be a r.v.\ with $\mathrm{GPD}(\xi,\psi)$
distribution in $(0,\sigma)$ and $c\geq0$ or $c\leq-\sigma$,
then $Y=-1/(X+c)$ has distribution $\mathrm{GPD}$ with location parameter
if and only if $c=\psi/\xi$. Then $Z=Y+1/c$ has $\mathrm{GPD}(-\xi,\xi^{2}/\psi)$
distribution.
\end{thm}

\begin{proof}
From \eqref{eq:gpd} the distribution function of $Y$ is
\begin{equation}
F_{Y}(y)\!=\!F(x(y))\!=\!1-\left(1-\frac{\xi}{\psi}\left(\frac{cy+1}{y}\right)\right)^{\!-1/\xi\!}\!=\!1-\left(\frac{\psi y}{y(\psi-\xi c)-\xi}\right)^{\!1/\xi},\label{eq:Fy}
\end{equation}
where $-1/c<y<-1/(\sigma+c)$. The denominator of the right
term of \eqref{eq:Fy} is a constant if and only if $c=\psi/\xi$.
In this case the distribution function of $Z$ is
\[
F_{Z}(z)=F_{Y}(y(z))=1-(1-\psi z/\xi)^{1/\xi}=1-(1-\xi z/(\xi^{2}/\psi))^{1/\xi},
\]
where $0<z<\sigma_{z}$, $\sigma_{z}=\xi/\psi$ for $\xi>0$ and $\sigma_{z}=\infty$
for $\xi<0$. Hence, $Z$~has $\mathrm{GPD}(-\xi,\xi^{2}/\psi)$
distribution and $Y$ has $\mathrm{GPD}$ distribution with location parameter.
\end{proof}

\begin{cor}
Let $\xi>0$, $\psi>0$ and $c=\psi/\xi$, then a r.v.\ $X$ has $\mathrm{GPD}(\xi,\psi)$
distribution if and only if $Z=X/(c(X+c))$ has $\mathrm{GPD}(\xi_{z},\psi_{z})$
distribution with $\xi_{z}=-\xi$, $\psi_{z}=\xi^{2}/\psi$ and the
support $(0,\xi/\psi)$.
\end{cor}

\begin{proof}
In the direct sense, this is proved by the Theorem~\ref{thm:9}, because
$c>0$ and $Z=X/(c(X+c))=-1/(X+c)+1/c$. 

The converse is also a consequence of Theorem~\ref{thm:9}, because
the inverse of the above transformation is 
\[
X=c^{2}Z/(1-cZ)=Z/(c_{2}(Z+c_{2}))=-1/(Z+c_{2})+1/c_{2}
\]
where $c_{2}=-1/c=-\xi/\psi$. The support of $Z$ is $(0,\psi_{z}/|\xi_{z}|)=(0,\xi/\psi)$
and $Z+c_{2}<0$ (equivalently $c_{2}\leq-\xi/\psi$), then $X$ is
a monotonous increasing function of $Z$ and Theorem~\ref{thm:9}
proves the result.
\end{proof}

In practical applications of the previous results, a first estimate of the shape and the scale parameters is required in order to define the transformation to a lighter tail, after which the residual empirical CV plot is constructed.

\section{Multiple threshold test}
\label{section3}

Some results of Castillo \emph{et al.}~\cite{Castillo14} on the residual CV extend directly from the exponential distribution to all GPD, provided there is a finite fourth moment. Therefore, the proof of the following theorem is omitted. The asymptotic distribution of the residual CV  as a random process indexed by the threshold provides pointwise error limits for $\mathrm{CV}$-plot in \eqref{eq:sample cv} and a multiple thresholds test for GPD that really does reduce the multiple testing problem. The multiple thresholds test provides a clear sense of significance levels and p-values.

\begin{thm}\label{thm:T1}
Let $X$ be a $\mathrm{GPD}(\xi,\psi)$ distributed
r.v., with $\xi\!<\!1/4$. Then $\sqrt{n}(cv_{n}(t)-c_{\xi})$, where $cv_{n}(t)$ and $c_{\xi}$ were respectively defined in  \eqref{eq:sample cv} and  \eqref{eq:cvgpd},
converges to a Gaussian process with zero mean and covariance function
given by
\[
\rho_{0}(s,t)=\exp((s\wedge t)/\psi),
\]
for $\xi=0$, and
\begin{equation*}
\begin{split}
\rho_{\xi}(s,t)&=(((\psi+\xi s)/\psi)^{1/\xi})(1-\xi)^{2}(6\xi^{4}t^{2}+12\psi\xi^{3}t+8\xi^{3}st
-9\xi^{3}t^{2}+6\psi^{2}\xi^{2}\\
&\quad+8\psi\xi^{2}s-10\psi\xi^{2}t
-2\xi^{2}st+3\xi^{2}t^{2}-\psi^{2}\xi-2\psi\xi s+4\psi\xi t+\psi^{2})\\
&\quad/((1-3\xi)(1-2\xi)^{2}(1-4\xi)(\psi+\xi s)^{2})
\end{split}
\end{equation*}
for $\xi\neq0$ and $s\leq t$.
\end{thm}

Pointwise error limits of the $\mathrm{CV}$-plot under GPD follow from
the next result.

\begin{cor}\label{cor:C4}
Given a sample $\{ X_{j}\} $ of a $\mathrm{GPD}(\xi,\psi)$
distribution $(\xi<0.25)$ and a threshold~$t$,  the asymptotic distribution of the residual $\mathrm{CV}$
is 
\begin{equation}
D_{t}(\xi)\equiv\sqrt{n(t)}(cv(t)-c_{\xi})\overset{d}{\rightarrow}N(0,\sigma_{\xi}^{2}).\label{eq:CI}
\end{equation}
where $c_{\xi}$ is in \eqref{eq:cvgpd},  $n(t)=\sum_{j=1}^{n}1_{(X_{j}>t)}$ and
\[
\sigma_{\xi}^{2}=\frac{(1-\xi)^{2}(6\xi^{2}-\xi+1)}{(1-2\xi)^{2}(1-3\xi)(1-4\xi)}.
\]
\end{cor}

\begin{proof}

The asymptotic variance is given by $\sigma_{\xi}^{2}=\rho_{\xi}(0,0)$,
where the covariance function is in Theorem~\ref{thm:T1}. The Theorem  can be applied to the threshold excesses $\{ X_{j}-t\mid X>t\} $, replacing $n$ with $n(t)$ and $cv(0)$ with $cv(t)$.
From \eqref{eq:invshape} the threshold excesses are again $\mathrm{GPD}$
with the same parameter~$\xi$ and the $\mathrm{CV}$ does not depend on~$\psi$.
\end{proof}

From the last result the asymptotic confidence intervals of the $\mathrm{CV}$-plot
for exponential distribution are obtained with $c_{0}=1$ and $\sigma_{0}^{2}=1$
and for uniform distribution with $c_{-1}=1/\sqrt{3}$ and $\sigma_{-1}^{2}=8/45$.

\subsection{ Simple null hypothesis}

Corollary~\ref{cor:C4} makes it possible to test
whether the empirical $\mathrm{CV}$ of a sample, or of threshold excesses, fit
the $\mathrm{CV}$ of a $\mathrm{GPD}$ with fixed values~$\xi$ and~$t$. 
However, from~\cite{Gupta00}, in order to have a consistent test in GPD, $CV(t)=c_\xi$  must be checked for all threshold  $t$. From  Theorem~\ref{thm:T1}, a multiple threshold test for a number $m$ of thresholds as large as necessary for practical applications can also be constructed using the building blocks 
$D_{t}^{2}(\xi)/\sigma_{\xi}^{2}$, regardless of the scale parameter, with asymptotic distribution $\chi_{1}^{2}$ under the null hypothesis
of~$\mathrm{GPD}$ $(\xi<0.25)$.

The choice  of thresholds could be arbitrary, but the multiple thresholds test, $T(\xi)$, is designed 
to avoid subjectivity as much as possible, to the limit of the number of thresholds $m$. 
If the thresholds are selected as empirical quantiles or order statistics, then  $T(\xi)$ is invariant when the sample is multiplied by a positive number while maintaining the set of probabilities, since CV is invariant. 
This first condition ensures that the test results do not depend on units used for the observations.
 
Given a sample $\{ x_{j}\} $
of size~$n$ of non-negative numbers, $Q_{n}(p)$ denotes
the inverse of the empirical distribution function,
\begin{equation}
Q_n(p)=\inf [x:F_n(x)\geq p].\label{eq:QuanFunc}
\end{equation}
From a set of probabilities $\{ 0=p_{0}<p_{1}<\dotsb<p_{m}\} $
let $\{ 0=q_{0}<q_{1}<\dotsb<q_{m}\} $ be the corresponding
empirical quantiles of the sample, $q_{k}=Q_{n}(p_{k})$,
then a multiple thresholds statistic can be constructed as
\[
T(\xi)=\sum_{k=0}^{m}D_{q_{k}}^{2}.
\]

The asymptotic expectation is $(m+1)\sigma_{\xi}^{2}$,
hence $T(\xi)/(m+1)$ is an estimator of the
asymptotic variance $\sigma_{\xi}^{2}$, when $\xi$ is known or estimated. Note that the distribution
of $T(\xi)$ is independent of the scale parameter~$\psi$.
$T(\xi)$ makes it possible to test the null hypothesis
that the sample comes from a distribution with the residual $\mathrm{CV}$ corresponding
to previous quantiles all equal to~$c_{\xi}$.
\[
H_{0}: \mathrm{CV}(q_{k})=c_{\xi},\quad k=0,1,\dotsc,m.
\]
Hence, if $H_{0}$ is accepted and $m$ is large enough, say $20$
or $50$, it will be reasonable to assume that the sample comes from
a distribution $\mathrm{GPD}(\xi,\psi)$. The previous test~$T(\xi)$
is a global test in the sense that some $D_{q_{k}}^{2}$ may be significant
and others not but with one test alone the equality of all $\mathrm{CV}$ for
all quantiles is checked. 

A second desirable condition is to select the set of probabilities that determine the statistic
$T(\xi)$ so that  the corresponding thresholds are approximately equally spaced. This can be achieved for the exponential distribution by taking
$0<p<1$, $p_{k}=1-p^{k}$, $(k=0,\dotsc,m)$ and $q_{k}$
the corresponding quantiles, since for a random variable, $X$, with exponential distribution
$\Pr\!\left\{ X\!>\left(\mu\log(1/p)\right)k\right\}\! =p^{k}$, where $\mu$ is the expected value. Then the condition holds for $\xi=0$ and is fairly approximate
for $\xi$. Selecting the probabilities this way, $q_{k}=Q_{n}(p_{k})\approx x_{(n-np^{k})}$,
$n(q_{k})\approx n\, p^{k}$ and $T(\xi)$
becomes

\begin{equation}
T_{m}(\xi)=n\sum_{k=0}^{m}p^{k}(cv(q_{k})-c_{\xi})^{2}.\label{eq:Tm}
\end{equation}

In applications, given the number of single tests that will be included
in the multivariant test, $m$, we choose the value of $p$, which
determines the distance between the quantiles, such that $n\, p^{m}\approx ns$,
where $ns$ is the sample size such that irrelevant information comes
from smaller sub-samples. Hence, given $m$, $p=(ns/n)^{1/m}$
is suggested. In this paper $ns\approx8$ is used in numerical algorithms.
Note that this way $T_{m}(\xi)$ depends only on~$\xi$
and $m$ and the researcher chose only the number of thresholds used
in the analysis, essentially eliminating subjectivity. These multiple
thresholds tests generalize those developed by Castillo \emph{et al.}~\cite{Castillo14} 
for $\xi=0$ and $p=1/2$.

The asymptotic distribution of $T_{m}(\xi)$ is easily
calculated from Theorem~\ref{thm:T1}, following the steps suggested
by Castillo \emph{et al.}~\cite{Castillo14}, 
whenever $\xi<0.25$. However, taking
into account the different values of the extreme value index
and the diverse small sample sizes, it is easier in practice to calculate
the $p$-value for $T_{m}(\xi)$ using simulation methods,
which are especially simple in this case. Assuming $\mathrm{GPD}$ for simulations,
only the sample size, the number of thresholds, $m$, and $\xi$ are
needed. Since the distribution does not depend on scale, parameter
$\psi=1$ will be used.

\subsection{ Composite null hypothesis}

In most cases the parameter~$\xi$ is unknown and its estimate should
be incorporated in the statistic~$T_{m}(\xi)$ (see the R
code below). The method for estimating~$\xi$ leads to
slight variations in the statistic, but it leads to essentially equivalent
inference whenever we use the same estimation method in simulations
to obtain the $p$-value. The null hypothesis is now that the sample
comes from a distribution in which all $(m+1)$ residual
$\mathrm{CV}$ are equal. 
\[
H_{0}: \mathrm{CV}(q_{0})=\dotsb=\mathrm{CV}(q_{m}),\quad k=0,1,\dotsc,m.
\]
The alternative hypothesis is that the residual $\mathrm{CV}$ are equal from
a threshold~$q_{r}$ $(0<r\leq m)$ to the threshold~$q_{m}$.

The most recommended estimation method is  maximum likelihood estimation (MLE),
although in $\mathrm{GPD}$ it is only asymptotically efficient provided $-0.5<\xi$,
see~\cite{Davison90}.  
For this distribution, the $\mathrm{CV}$ is a one-to-one
transformation of $\xi$, see \eqref{eq:cvgpd}, and the empirical
$\mathrm{CV}$ of the residual sample, $\mathrm{CV}(t)$, provides an alternative method of estimation.
It is asymptotically normal whenever $\xi<0.25$, see Corollary~\ref{cor:C4}.
The multiple thresholds tests \eqref{eq:Tm} suggest estimating $\xi$
as the value such that $c_{\xi}$~achieves the minimum $T_{m}(\xi)$,
namely

\begin{equation}
\tilde{c}_{\xi}=\sum_{k=0}^{m}p^{k}cv(q_{k})/\sum_{k=0}^{m}p^{k}=(1-p)\sum_{k=0}^{m} p^{k}cv(q_{k})/(1-p^{m+1}).\label{eq:estim}
\end{equation}

From Corollary~\ref{cor:C4} the estimator is also asymptotically normal. The main
advantage of this method is that under the alternative hypothesis
it is a better estimator than $\mathrm{CV}$ or MLE, since the sample is only $\mathrm{GPD}$ over a threshold~$q_{r}$. 
Since the main interest is in samples that are not $\mathrm{GPD}$, but in the
tail, and results are often used in small samples with $\xi<0$, the estimation method \eqref{eq:estim} is included in the statistic $T_m=T_{m}(\tilde{\xi})$. The following
R code for~$T_{m}$ is used in the algorithms, see~\cite{RCore10}.  

\begin{verbatim}
#Statistic Tm of a sample given the number of thresholds m.
Tm<-function(m,sample){sam<-sample-min(sample);
     n<-length(sam);ns<-8;
     p<-round(exp(log(ns/n)/m),digits=2);
     Ws<-Ps<-Qs<-Cs<-numeric(m+1);
     for(k in 1:(m+1)){Ws[k]<-p^(k-1)};
     Ps<-1-Ws;Qs<-as.vector(quantile(sam,Ps));
     for(k in 1:(m+1))
     {Cs[k]<-sd(sam[sam>=Qs[k]]-Qs[k])/mean(sam[sam>=Qs[k]]-Qs[k])};
     cx<-(1-p)*sum(Ws*Cs)/(1-p^(m+1));xi<-(cx^2-1)/(2*cx^2);
     tm<-n*sum(Ws*(Cs-cx)^2);list(CV=cx,Tm=tm,Xi=xi)}
\end{verbatim}

\subsection{Threshold Selection Algorithms}\label{section4}

To select the number of extremes used in applying the peaks over a high threshold method,
threshold selection algorithms are developed in this section to estimate the point above which the 
$\mathrm{GPD}$ distribution can be used to estimate the extreme value index for a set of extreme events, $\{ x_{j}\} $, of size~$n$. For this purpose the previous statistical tests will
be adapted.

Note that in the $T_m$ calculation the number of thresholds~$m$ is
the only parameter that must be fixed by the researcher. This determines
the thresholds (quantiles) where the $\mathrm{CV}$ is calculated, $\{ q_{0}<q_{1}<\dotsb<q_{m}\} $,
which are fixed throughout the procedure. Then, by simulation of $\mathrm{GPD}$,
the associated $p$-value is calculated (running $10^{4}$~samples).
After that, we accept or reject the null hypothesis with the estimated
shape parameter using all the thresholds. 

If the hypothesis is rejected, the threshold excesses $\{ x_{j}-q_{1}\} $
are calculated for the sub-sample $\{ x_{j}\geq q_{1}\} $.
The previous steps are repeated, but removing one threshold, to accept
or reject the null hypothesis that the sample is from a $\mathrm{GPD}$. At every
stage only statistics associated to thresholds $k=r,\dotsc,m$, where
$0\leq r\leq m$, are calculated:
\begin{equation}
T_{m}^{r}(\xi)=n\sum_{k=r}^{m}p^{k}(cv(q_{k})-c_{\xi})^{2}.\label{eq:Tm-1}
\end{equation}

In summary, the \textbf{steps of the general algorithm} are
\begin{enumerate}
\item[(1)] Given $m$ find $p$ such that $np^{m}\approx ns$, where $ns$ is
the smaller sample size used to calculate $\mathrm{CV}$ (here $ns=8$ is used,
but it can be modified).
\item[(2)] Calculate $\{ 0=p_{0}<p_{1}<\dotsb<p_{m}\} $, where $p_{k}=1-p^{k}$,
and $\{ 0=q_{0}<q_{1}<\dotsb<q_{m}\} $, where $q_{k}=Q_{n}(1-p^{k})$,
$k=1,\dotsc,m$.
\item[(3)] Estimate $\tilde{\xi}$ minimizing the value of $T_{m}(\xi)$
with the specific values in the previous steps.
\item[(4)] Calculate by simulation of $\mathrm{GPD}$ the $p$-value associated to the min\-i-\linebreak mum~$T_{m}(\tilde{\xi})$ and accept or reject the null
hypothesis with the estimated shape parameter using all the thresholds
(starting with $q_{0}=0$). 
\item[(5)] If the hypothesis is rejected, compute the threshold excesses $\{ x_{j}-q_{1}\} $
for the sub-sample $\{ x_{j}\geq q_{1}\} $ and repeat
the previous steps with $\{ p_{1}<\dotsb<p_{m}\} $ and $\{ q_{1}<\dotsb<q_{m}\} $,
to accept or reject the null hypothesis that the sample is from a
$\mathrm{GPD}$, but removing a threshold.
\item[(6)] Continue the process for the next value in the index of thresholds
while the hypothesis is rejected.
\end{enumerate}
Several authors recommend giving a prominent role to the exponential
distribution in the model $\mathrm{GPD}$, see \cite{CastilloSerra14}.  
The usual
method for doing this is to consider the exponential models as the
null hypothesis testing against $\mathrm{GPD}$, see \cite{Kozubowski09}.  
Alternatively, one can consider the Akaike or Bayesian information
criteria for model selection, see \cite{Clauset09}. 
The previous algorithm can be adapted to the case when $\xi=0$ (or simply
known) skipping step-3.

\section{Danish fire insurance data}\label{exa:Danish}

An interesting aspect of this article is the combination of the results of sections ~\ref{section2} and ~\ref{section3} when applying the  peaks over threshold technique for tails in any maximum domain of attraction. This approach is illustrated here using a popular dataset.

The Danish fire insurance data are a well-studied
set of losses to illustrate the basic ideas of extreme value theory.
The dataset consists of $2,156$~fire insurance losses over one million
Danish kroner from~1980 to~1990 inclusive, see 
\cite[Example~6.2.9]{Embrechts97}, \cite{Resnick97} 
and  \cite[Example~7.23]{McNeil05}.

In this example the authors agree to assume iid observations and a
heavy tailed model. They also agree to set the threshold at $u=10$ million Danish kroner, the exceedances over the threshold, denoted $\{ x_{j}\} $, are $n_{10}=109$.

Fitting a $\mathrm{GPD}$ to $\{ x_{j}\} $ by MLE, the parameter estimates
in~\cite{McNeil05} 
are $\hat{\xi}=0.50$ and $\hat{\psi}=7.0$
with standard errors~$0.14$ and $1.1$, respectively. Thus the fitted
model is a very heavy-tailed, infinite-variance model and the method
in Section \ref{section3} cannot be applied directly. However, they can
be used through the results shown in Section~\ref{section2}. 

First of all, let us suppose we want to use $\mathrm{CV}$ to check whether the
above data correspond to a $\mathrm{GPD}$ distribution with the estimated extreme value index.
Applying Theorem~\ref{thm:9} with $c=\hat{\psi}/\hat{\xi}=14$,
let $z_{j}=-1/(x_{j}+c)+1/c$ be, then the set $\{ z_{j}\} $
has light tails and the same extreme value index with the sign changed, provided
that the estimated parameters are the true parameters. The $\mathrm{CV}$ of $\{ z_{j}\} $
is $cv=0.697$ which provides a new estimation of $\xi$, solving
\eqref{eq:cvgpd} by $\xi_{z}=(cv^{2}-1)/(2cv^{2})=-0.530,$
then, according to Theorem~\ref{thm:9}, $\tilde{\xi}=-\xi_{z}=0.53$,
not far from $0.50$, since the standard error is~$0.14$. Alternatively,
the multiple thresholds statistic~$T_m$, from \eqref{eq:Tm-1},
can be used to check $\xi=0.5$. The corresponding $\mathrm{CV}$ under $\mathrm{GPD}$ is
$c_{\xi}=0.707$. Taking $m=20$, we get $T_m=4.89$ with a $p$-value
$0.421$ (by simulation with $10^{4}$ samples), accepting the null
hypothesis.

Now consider the problem of choosing the threshold to estimate the
extreme value index. In this example, most researchers use a visual observation
of the \emph{ME}-plot on the full Danish dataset. The algorithm in
Section~\ref{section4} with the transformations from Section~\ref{sec:Transformations}, comes to similar
solutions automatically and opens up new perspectives.

\begin{figure}[ht]
\centering
\begin{tabular}{|c|c|}
\hline 
\includegraphics[scale=0.32]{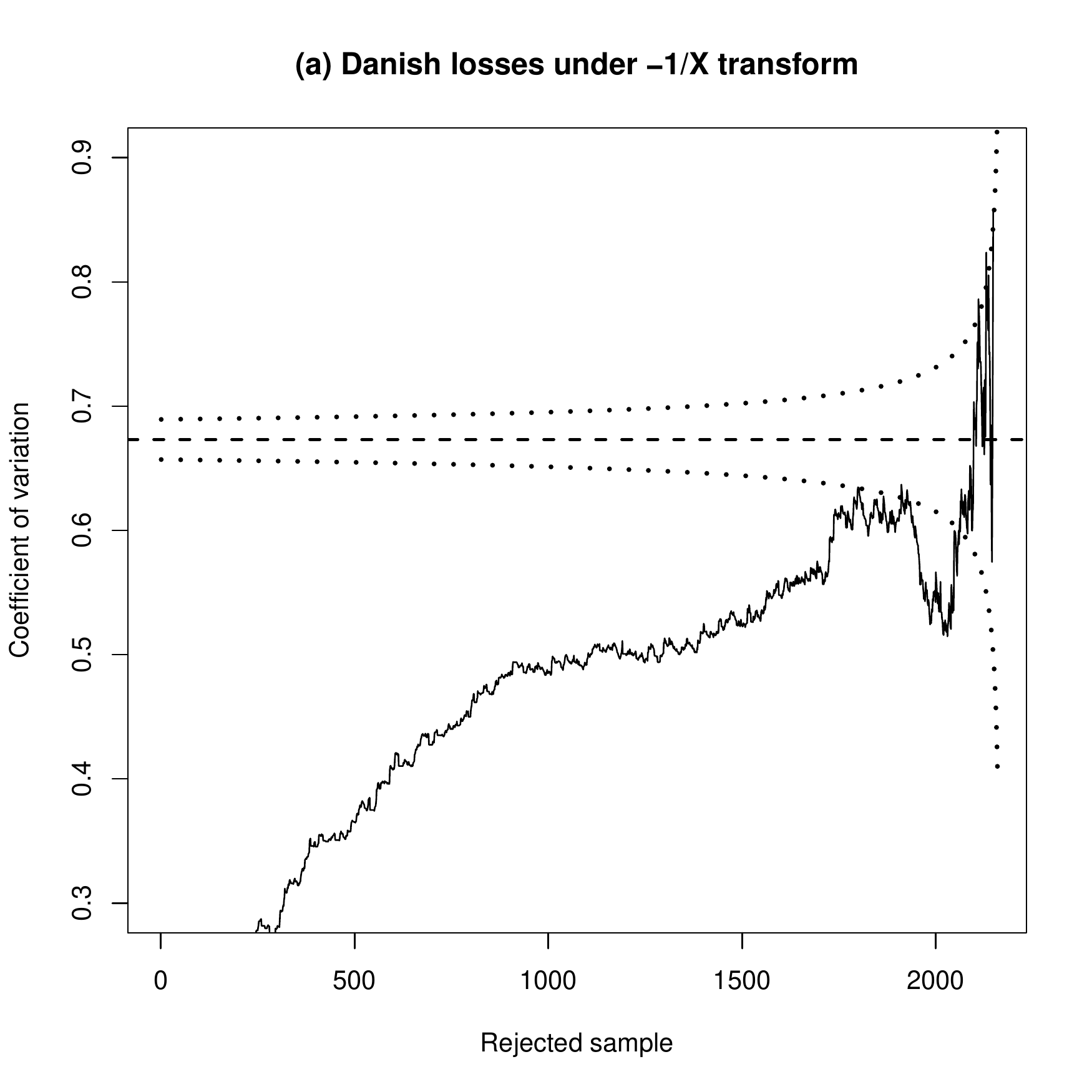} & \includegraphics[scale=0.32]{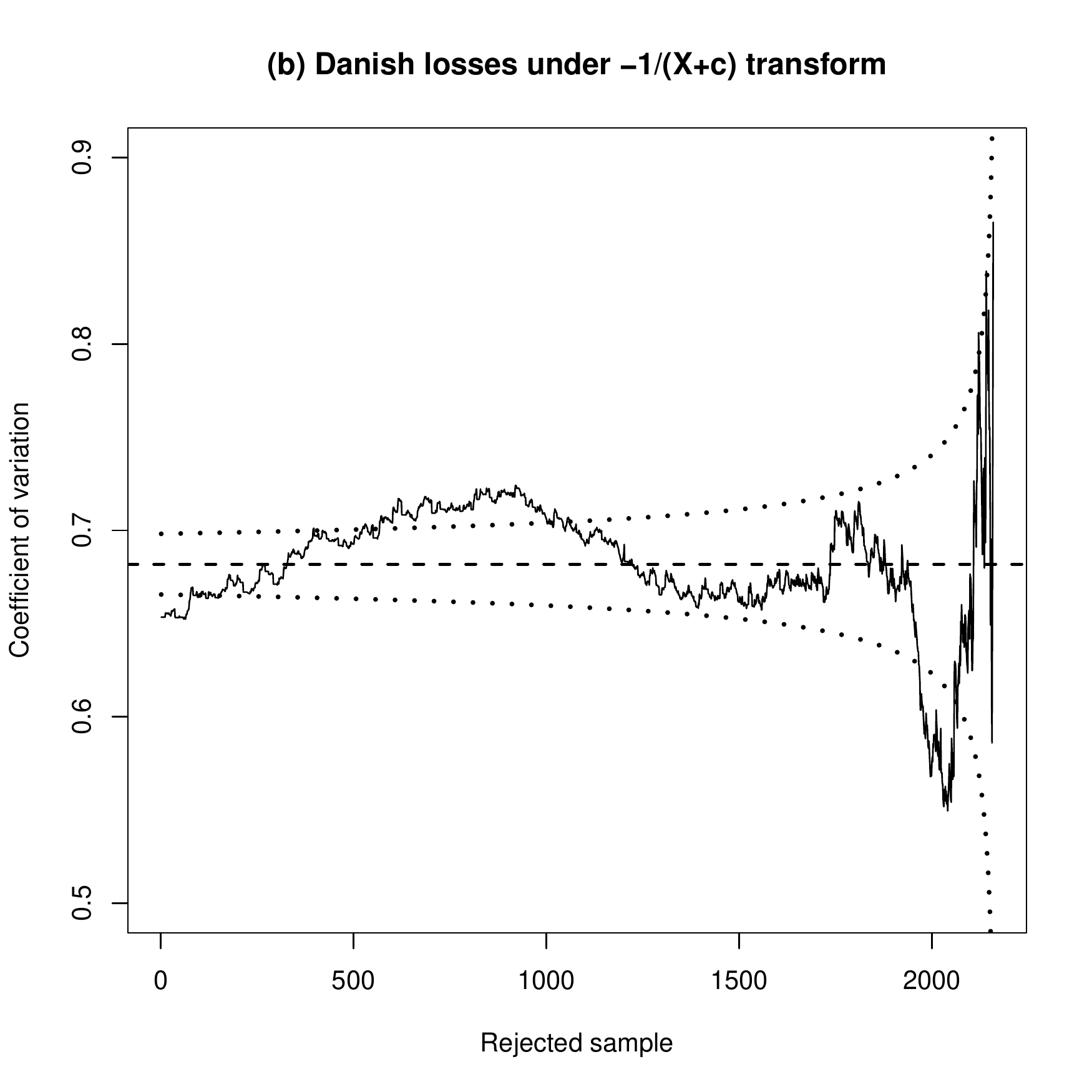}\\
\hline 
\end{tabular}
\caption{Residual empirical $\mathrm{CV}$ for The
Danish fire insurance losses under transformation of the data. (a):~Dataset, transformed by $-1/X$ . (b):~Dataset, transformed by $-1/(X+\psi/\xi)$.
The dotted lines correspond to the asymptotic confidence intervals
($90\%$) under the estimated parameter, the dashed line is its $\mathrm{CV}$.}
\label{fig:F5}
\end{figure}

Figure \ref{fig:F5} shows the $\mathrm{CV}$-plots of the full Danish
dataset, transformed according to the Corollary~\ref{cor:8}, plot~(a), and Theorem~\ref{thm:9}, plot~(b). The first, corresponding to the transformation $X^{*}=-1/X$, shows an increasing $\mathrm{CV}$ and the
second, corresponding to $Z=-1/(X+c)+1/c$, shows a stabilized $\mathrm{CV}$
close to a constant, indicating that the original dataset is close
to a $\mathrm{GPD}$, which is also shown by \emph{ME}-plot.

Applying the algorithm of Section~\ref{section4} with $m=20$ after transformation~$X^{*}$, constant residual $\mathrm{CV}$ is rejected in the first 11 steps (each
one reduces the sample size by $(1-p)=24\%$). Step 12,
for the last 106 observations, accepts constant residual $\mathrm{CV}$ ($p$-value
$= 0.269$) with estimates $c_{\xi}=0.673$ and $\xi=0.603$. The
estimated threshold is approximately the same ($u=10.2$ instead of~$10$), while the extreme value index is different but within the confidence interval.

The algorithm in Section~\ref{section4}, with $m=20$ after transformation $Z$
with $c=0.932/0.611=1.524$, rejects constant residual $\mathrm{CV}$ in the first
three steps. Step 4, for the last 951 observations, accepts constant
residual $\mathrm{CV}$ ($p$-value $= 0.167$) with estimates $c_{\xi}=0.675$
and $\xi=0.599$. The number of observations is much higher, the extreme value index
being very close to that obtained with the transformation~$X^{*}$
and within the confidence interval. The $p$-value remains similar
in the following steps up until the 12th, where it jumps up to $0.474$.
The number of observations is again~$106$ and the estimation $\xi=0.548$,
nearer to~$0.50$.

The conclusions from using the new methodology to analyze this dataset
are the following. First, the results obtained by previous investigators
are validated, in particular $\mathrm{GPD}$ can be accepted with parameter $\xi=0.5$,
for the $109$~larger observations see \cite{McNeil05}. 
This also shows the consistency of the presented methodology with
other common techniques.

Moreover, from examining the extreme value index it is now known that for the
$951$~larger observations $\mathrm{GPD}$ can also be accepted, where the MLE
parameter estimate is $\xi=0.680$, with standard error $0.055$ ($\xi=0.599$
obtained by $T_m$ is within the confidence interval). The estimated
extreme value index is now much more accurate because the sample size is much
larger. We also note that the tails are heavier than was assumed,
which means that higher risks should be considered.

\section*{Acknowledgement}

This work was supported by the Spanish Ministry of Economy and Competitiveness
under Grant: Applied Stochastic Processes, MTM 2012-31118.

\end{document}